\documentclass[reqno]{amsart}

\usepackage{amsmath,amssymb,amsthm}

\usepackage{tikz}

\usepackage{caption}
\usepackage{url}
\usepackage{graphicx}
\usepackage{graphics}
\usepackage{algorithm}
\usepackage{algorithmicx}
\usepackage{algpseudocode}
\usepackage[backend=biber, url=false, doi=false, maxnames=10, maxalphanames=5]{biblatex}
\DeclareMathOperator*{\smax}{smax}

\begin{document}
\title[]{Low-Discrepancy Set Post-Processing\\ via Gradient Descent}

\author[]{Fran\c{c}ois Cl\'ement, Linhang Huang, Woorim Lee,\\ Cole Smidt, Braeden Sodt \and Xuan Zhang}

\address{Department of Mathematics, University of Washington, Seattle}
 \email{fclement@uw.edu }
 \email{lhhuang@uw.edu}
 \email{wlee1130@uw.edu}
\email{csmidt@uw.edu}
\email{bsodt@uw.edu}
\email{xuanz24@uw.edu}

 \begin{abstract}
     The construction of low-discrepancy sets, used for uniform sampling and numerical integration, has recently seen great improvements based on optimization and machine learning techniques. However, these methods are computationally expensive, often requiring days of computation or access to GPU clusters. We show that simple gradient descent-based techniques allow for comparable results when starting with a reasonably uniform point set. Not only is this method much more efficient and accessible, but it can be applied as post-processing to any low-discrepancy set generation method for a variety of standard discrepancy measures. 
 \end{abstract}
 \maketitle
\section{Introduction}

Discrepancy measures are a family of functions that quantify how well a given point set in $[0,1)^d$ approximates the uniform distribution. For all these measures, the closer the discrepancy gets to 0, the more uniformly distributed the point set is. Finite point sets or infinite sequences which minimize the discrepancy are called low-discrepancy sets and sequences. Research has focused on two aspects: understanding what is the minimal error made when approximating the uniform distribution, and constructing low-discrepancy sets that allow for efficient sampling techniques. On the theoretical side, despite extensive work~\cite{Chaz,Pano,Mat,Nie92}, only dimension 1 and dimension 2 for sets have been solved up to a constant factor, with a minimal error of $\mathcal{O}(\log(n)/n)$, which is reached by constructions such as the van der Corput and Sobol' sequences~\cite{vdC,Sobol}. In higher dimensions $d$, the optimal order is still a wide open question, with a lower bound of $\Omega(\log(n)^{c(d)+(d-1)/2})$ for a point set of size $n$ obtained by Bilyk, Lacey and Vagharshakyan~\cite{BilykSmall}. This should be compared to the best current upper-bound obtained by many known constructions, $\mathcal{O}(\log^{d-1}(n)/n)$.

On the practical side, low-discrepancy constructions have been ubiquitous in numerical integration since the results of Koksma and Hlawka~\cite{Hlawka,Koksma}, who showed that the error made when approximating the integral of a function $f$ with a finite sum of local evaluations taken in some points $P$ depends directly on a property of $f$ and on the $L_{\infty}$ star discrepancy of the point set $P$. Since $f$ can only very rarely be modified, decreasing the discrepancy of $P$ is the simplest solution to improve the numerical approximation of integrals. Apart from numerical integration, low-discrepancy constructions have also seen widespread use in computer vision~\cite{MatBuilder}, financial mathematics~\cite{GalFin}, experimental design~\cite{SantnerDoE}, optimization~\cite{CauwetCDLRRTTU20,DiedOpti} or motion planning~\cite{RUSCH2024}.

The $L_{\infty}$ star discrepancy is arguably the most prominent discrepancy measure, largely due to its role in the Koksma-Hlawka inequality. For a fixed point set $P$ in $[0,1]^d$, it consists in finding the worst absolute difference between the number of points of $P$ inside a box anchored in $(0,\ldots,0)$ and the volume of this box. More formally, it is given by
$$d^*_{\infty}(P):=\sup_{q\in [0,1)^d}\left|\frac{|P\cap [0,q)|}{|P|}-\lambda([0,q))\right|,$$
where $\lambda$ is the Lebesgue measure of $[0,q)$. Despite its theoretical usefulness, the $L_{\infty}$ star discrepancy is unwieldy to use in practice: simply evaluating it is known to be NP-hard~\cite{complexity}, and even $W[1]$-hard in the dimension~\cite{W1hard}. The best algorithm known to this has a $\mathcal{O}(n^{1+d/2})$ runtime~\cite{DEM}, while a threshold accepting-based heuristic allows evaluation for dimensions higher than 10~\cite{GnewuchWW12}. Unsurprisingly, optimizing the \emph{construction} of sets with low $L_{\infty}$ star discrepancy is an even harder problem. Until recently, constructions relied either on purely mathematical approaches (see~\cite{Nie92} for traditional constructions and bounds), or were limited to \emph{extremely} small set sizes~\cite{Larcher, PTV,whit:onop:1976}.

Recent computational methods to optimize the construction have globally relied on two approaches. The first has been to directly optimize the $L_{\infty}$ star discrepancy via non-linear programming~\cite{PnasCDKP,CDKP} or algorithmic techniques and heuristics~\cite{CDP23,CDP}. These methods tend to be computationally very expensive, and are quickly limited by the cost of evaluating the $L_{\infty}$ star discrepancy. The second has been to leverage another discrepancy measure, the $L_2$ star discrepancy (introduced in more detail in Section 2). This second measure is \emph{much} easier to evaluate and, for reasons not perfectly well understood yet, seems to approximate extremely well the $L_{\infty}$ star discrepancy in low dimensions. The Message-Passing Monte Carlo method in~\cite{MPMC} was the first to really leverage this observation, before being followed by~\cite{SubsetL2}. These methods do remain computationally expensive, requiring high computing power for~\cite{MPMC} and potentially multiple hours for a single run of~\cite{SubsetL2}. With the aim of constructing infinite sequences, and not finite sets, the $L_2$ star discrepancy has also led to surprisingly good results in one dimension~\cite{CKritz,Kritz}.

\begin{figure}[t!]
\includegraphics[width = 0.7\linewidth]{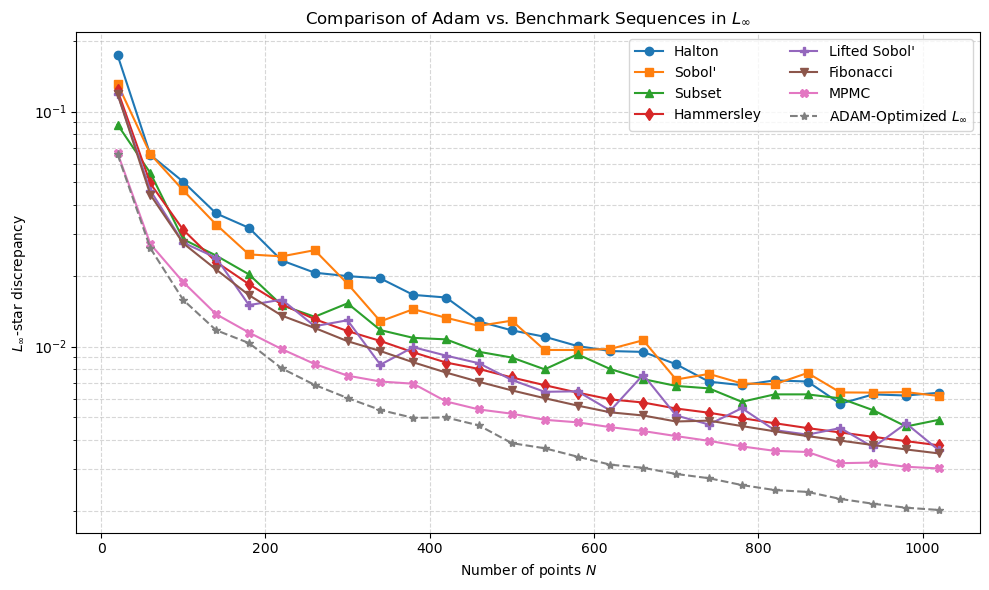}
\caption{The $L_\infty$ values for $n=20$ up to $n=1020$ for well-known low discrepancy sets. Projected Gradient Descent is in gray.}
\label{fig:Linftycomp}
\end{figure}

\textbf{Our Contribution.} In this paper, we aim to show that much simpler tools can be used to optimize the $L_2$ star discrepancy and leverage its connection to the $L_{\infty}$ star discrepancy. We implement an ADAM-based projected gradient descent algorithm for the $L_2$ star discrepancy as well as other common $L_2$ discrepancies such as the periodic discrepancy. While the gradient descent can not be applied to \emph{any} point set, we show that starting with classical sets such as the Fibonacci lattice lead to low-discrepancy sets comparable or better than state-of-the-art results, even for the $L_{\infty}$ star discrepancy. Figure~\ref{fig:Linftycomp} shows that for the $L_{\infty}$ star discrepancy in two dimensions, the gradient descent algorithm outperforms the results from~\cite{MPMC}. Furthermore, our methods can be applied to any set obtained from more sophisticated methods such as~\cite{SubsetL2,CDP23,CDP,MPMC}. Empirically, we observe that even when the initial set is near optimal, the gradient descent algorithm never makes the point set worse, and usually manages to improve it at least marginally. The very low runtime, between a few seconds and a minute on a laptop, as well as its flexibility regarding different discrepancy measures, make this method a very interesting tool for post-processing existing sets. All our code is available at \url{https://github.com/linhang-h/L2_Discrepancies}.

\textbf{Structure of the paper.}
Section 2 describes the $L_2$ star discrepancy we are optimizing, as well as a wider class of discrepancy functions which have a similar structure and can be optimized in a similar way. Section 3 presents our algorithm and the modifications that were required to adapt gradient descent to $L_2$ discrepancies. Section 4 describes which input parameters (initial point set and number of steps) lead to good algorithm performance. Section 5 then presents our results, including a comparison with previous methods and some results for other discrepancy measures we hope will serve as benchmarks for future improvement.

\section{$L_2$ discrepancies}
\subsection{The $L_2$ star discrepancy}

Like the $L_{\infty}$ star discrepancy, the $L_2$ star discrepancy only considers boxes anchored in the origin. However, rather than considering the absolute worst difference, it describes the average absolute difference. More precisely, for a point set $P:=(x_i)_{i \in \{1,\ldots,n\}}$ in dimension $d$,

$$d_{2}^*(P):=\int_{q \in [0,1)^d} \left|\frac{|P\cap [0,q)|}{|P|}-\lambda([0,q))\right|dq.$$
While this formula does not seem much easier to work with at a glance, a result of Warnock~\cite{Warnock} shows that it can be rewritten into
$$d_2^*(P)=\frac{1}{3^d}-\frac{2^{1-d}}{n}\sum_{i=1}\prod_{k=1}^d(1-x_{i,k}^2)+\sum_{i,j=1}^n\prod_{k=1}^d(1-\max(x_{i,k},x_{j,k})),$$
where $x_{i,k}$ is the $k$-th coordinate of $x_i$. This expression has the advantages of both clearly identifying the contribution of each point and, apart from the $\max$ terms, of being differentiable. However, as we will describe in Section 3, it is well-known how to handle such terms to be able to apply gradient descent techniques. 
It is also much faster to compute the $L_2$ star discrepancy than the $L_{\infty}$ star discrepancy, with the Warnock formula requiring only $\mathcal{O}(dn^2)$ operations. We will not be using it here, but we note that there exists another algorithm by Heinrich that requires $\mathcal{O}(n\log(n)^d)$ operations~\cite{HeinL2}.

\subsection{Other $L_2$ discrepancies}

While star discrepancies are generally the main area of focus due to the Koksma-Hlawka inequalities and similar inequalities derived from it~\cite{DickP10,HickQuad}, they have some known drawbacks. The main one for applications is the bias induced by the anchoring: the corner $(0,\ldots,0)$ of the box $[0,1)^d$ becomes more important than the others. However, for applications where the parameter space is rescaled to obtain $[0,1)^d$, it usually does not make sense to make the minimal parameter values more important. Other measures such as the extreme discrepancy (all axis-parallel boxes, not necessarily anchored in $(0,\ldots0)$), or the periodic discrepancy (all axis-parallel boxes where $[0,1]^d$ is associated with the $d$-dimensional torus) have been regularly studied theoretically~\cite{KriHinPill,Mat}, but dedicated constructions are rare. 

In the $L_{\infty}$ setting, they are even more computationally expensive than the star discrepancy and, to our knowledge, there is no dedicated algorithm to compute them (bar repeated use of the algorithm for the star discrepancy in~\cite{DEM}, or applying the structural observation in~\cite{NieBox}). In the $L_2$ setting however, all these other discrepancy measures also have a Warnock-like formula, as summarized recently in~\cite{clément2025optimizationdiscrepancymeasures}. For example, a result of~\cite{KriHinPill} shows that the $L_2$ periodic discrepancy is given by 

$$d_{2}^{per}(P):=-\frac{1}{3^d}+\frac{1}{n^2}\sum_{i,j=1}^n\prod_{k=1}^d \left(\frac{1}{2}-|x_{i,k}-x_{j,k}|+(x_{i,k}-x_{j,k})^2\right).$$

Only the absolute value needs to be modified in this formula for us to be able to apply gradient descent: just like~\cite{clément2025optimizationdiscrepancymeasures} adapts the methods from~\cite{MPMC} to a wider range of discrepancy measures, we can adapt the gradient descent algorithm for $L_2$ star without radically changing the framework.

\section{Description of the algorithm}

This section describes the gradient descent algorithm and implementation for the $L_2$ star discrepancy. It can directly be modified for other discrepancy values by replacing the Warnock formula with the appropriate discrepancy function, and using a well-chosen smoothing function for the non-differentiable part of the formula.

In the Warnock formula, the only non-differentiable part is the maximum function that appears in each of the pairwise contributions. We replace the maximum by the $\tau$-softmax, given by
$$\smax(a,b,\tau):=\frac{1}{2}\left(a+b+\sqrt{(a-b)^2+\tau}\right).$$

In our experiments, we take $\tau=10^{-15}$. Given a set of points $P=(x_i)_{i \in \{1,\ldots,n\}}$, the new formula for the loss function of our gradient descent algorithm is now given by
$$d_2^*(P,\tau):=\frac{1}{3}^d-\frac{2^{1-d}}{n}\sum_{i=1}^n\prod_{k=1}^d(1-x_{i,k}^2)+\frac{1}{n^2}\sum_{i,j=1}^n\prod_{k=1}^d(1-\smax(x_{i,k},x_{j,k},\tau).$$

One other element to take into account is that we require all our points to stay in $[0,1]^d$. For this, we use projected gradient descent where, if any point is pushed outside the box for one or multiple coordinates, these coordinates are set to $0$ or $1$ (depending on which is the nearest). Finally, the stepsize could vary greatly based on $n$ and $d$. To avoid having to find the appropriate one and how it should evolve during the algorithm, we use the ADAM optimizer~\cite{ADAM}. This optimizer uses past gradient information to automatically adjust the stepsize for each coordinate of each point. Using existing implementations in Python also has the added benefit of automatizing the gradient computation. Note that, naïvely, the gradient can be computed in $\mathcal{O}(d^2n^2)$ time, since for each coordinate of each point there will be $n+1$ different products to compute.

We now describe the main process of the algorithm, described in Algorithm~\ref{algo}. At each step $t$ of the gradient descent, we have a point set $P_t$. We compute its gradient $$ g_{t}=\nabla_{P}d_{2}^{\ast}(P_{t};\tau),$$
as well as its moments
\begin{align*}
  m_{t} &= \beta_{1} m_{t-1} + (1-\beta_{1})\,g_{t},\\
  v_{t} &= \beta_{2} v_{t-1} + (1-\beta_{2})\,g_{t}^{\odot 2}.
\end{align*}

Here, \(m_{t}\) tracks the exponential moving average of the gradient, while \(v_{t}\) tracks that of the element-wise square. $\beta_{1}$ and $\beta_2$ are parameters that are automatically optimized by the ADAM optimizer. The symbol ``\(\odot\)’’ denotes the entry-wise product. To avoid bias during the initial iterations, these moments are rescaled
$$\hat m_{t} = \frac{m_{t}}{1-\beta_{1}^{t}},\qquad
  \hat v_{t} = \frac{v_{t}}{1-\beta_{2}^{t}},$$
which finally leads to the expression of the the new point set pre-projection on $[0,1]^d$
$$P_{t+\frac{1}{2}}:=P_{t}-\alpha\frac{\hat m_{t}}{\sqrt{\hat v_{t}}+\varepsilon}.$$ 
The last step is to project these points on $[0,1]^d$ when necessary: for each $x'_{i,k} \in P_{t+\frac{1}{2}}$, $x_{i,k} \in P_{t+1}$ is given by $x_{i,k}=\min\left(1,\max(0,x'_{i,k})\right)$. We abbreviate this process as
$$  P_{t+1} = \Pi_{[0,1]^{2}}\!\bigl(P_{t+\tfrac12}\bigr)$$
in the pseudo-code.

\begin{algorithm}[h]
\caption{Projected gradient descent algorithm for the $L_2$ star discrepancy}
\label{algo}
\begin{algorithmic}[1]
\Require initial lattice \(P_{0}=P_{n}\), steps \(T\), learning rate \(\alpha\)
\State \(m_{0}\gets0,\;v_{0}\gets0\)
\For{$t=0$ to $T-1$}
  \State \(g_{t}\gets\nabla_{P} d_{n,2}^{\ast}(P_{t};\tau)\)
  \State Set $m_{t} = \beta_{1} m_{t-1} + (1-\beta_{1})\,g_{t}$ and $\hat m_{t} = \frac{m_{t}}{1-\beta_{1}^{t}}$.
  \State Set $v_{t} = \beta_{2} v_{t-1} + (1-\beta_{2})\,g_{t}^{\odot 2}$ and $\hat v_{t} = \frac{v_{t}}{1-\beta_{2}^{t}}$.
  \State Compute $P_{t+\tfrac12}=P_{t}-\alpha\frac{\hat m_{t}}{\sqrt{\hat v_{t}}+\varepsilon}$
  \State \(P_{t+1}\gets\Pi_{[0,1]^{2}}\!\bigl(P_{t+\tfrac12}\bigr)\)
\EndFor
\State \Return \(P_{T}\)
\end{algorithmic}
\end{algorithm}

\section{Correctly Initializing the Algorithm}

For the experimental results described in the following sections, experiments were run on a laptop with AMD Ryzen 7 5800H CPU (8 cores, 16 threads, 3.2 GHz) and 16 GB of RAM, running Windows 11 (64-bit). All our code was implemented in Python, with $L_{\infty}$ star discrepancy computations in higher dimensions in C using code from~\cite{DiscGECCO} inspired from code by Magnus Wahlström. For 200 iterations, each experiment took between a few seconds for $n\sim 100$, to around a minute for $n \sim 2000$. This should be compared with the hour(s) potentially necessary with subset selection~\cite{SubsetL2} or optimization techniques~\cite{PnasCDKP}, or the clusters necessary to properly optimize in~\cite{MPMC}. Finally, the Sobol' sequence was generated using the GNU Scientific Library.

\subsection{Initial Point Set.} Gradient descent naturally requires an initial point set. It turns out that this is the most challenging part of the algorithm: the quality of the output very strongly depends on the quality of the input. We tested our algorithm in two dimensions on traditional low-discrepancy point sets such as prefixes of the Sobol' sequence~\cite{Sobol}, the lattice set $$\left\{(i/n,\{\sqrt{2}i\}): i\in\{0,\ldots,n-1\}\right\},$$ and the Fibonacci point set defined by
$$\left\{(i/n,\{\varphi i\}): i \in \{0,\ldots,n-1\}\right\},$$
where $\varphi$ is the golden ratio and $\{x\}$ is the fractional part of $x$. The Fibonacci set in particular is known for its incredible regularity and the one-dimensional sequence associated with it is the best performing traditional construction for the $L_{\infty}$ star discrepancy. We also generated random point sets to compare.

\begin{table}[h! tbp]
\begin{center}
\begin{tabular}{| c | c | c | c | c |}
\hline
Initialization & $d_2^*$ & $d_2^*$  After PGD & $d_\infty^*$ & $d_\infty^*$  After PGD \\
\hline
Fibonacci & 0.003438 & 0.001893 & 0.01200 & 0.007035\\
$\sqrt{2}$ Lattice & 0.003714 & 0.001927 & 0.01192 & 0.007960\\
Sobol & 0.003525 & 0.002344 & 0.01546 & 0.01015\\
\hline
\end{tabular}
\caption{Star discrepancy values before and after Projected Gradient Descent (PGD) for several initialization sets, with $n=260$.}
\label{tab:1}
\end{center}
\end{table}

Table~\ref{tab:1} describes the results obtained for the deterministic sets with $n=260$. There is a clear improvement for all sets, in particular for the lattice-like constructions. Despite starting with a similar $L_2$ discrepancy value, the Sobol' sequence leads to poorer values for the $L_2$ and $L_{\infty}$ star discrepancies. Table~\ref{tab:2} describes the results obtained for random point sets. We note that for $2\,000$ random initial sets, the mean $L_2$ star discrepancy was $0.007789$, and the minimum 0.004529. Again, while there is a clear improvement, we are quite far from the results of the Fibonacci set or even those of the Sobol' sequence.

\begin{table}[h! tbp]
\begin{center}
\begin{tabular}{| c | c | c | c | c |}
\hline
Discrepancy & Median & Mean & Min & Max \\
\hline
$L_2^*$ & 0.003642 & 0.003618 & 0.003041 & 0.004411 \\
$L_\infty^*$ & 0.02238 & 0.01943 & 0.01432 & 0.04615\\
\hline
\end{tabular}
\caption{Results from a test of 200 initialization sets with uniformly chosen point coordinates.}
\label{tab:2}
\end{center}
\end{table}

We conjecture that these different results come from the algorithm struggling to change the relative positions of the points. If we have points $x_{i}$ and $x_{j}$ such that for the $k$-th coordinate $x_{i,k}\leq x_{j,k}$, then at the end of the algorithm we will \emph{usually} have $x_{i,k}^{end}\leq x_{j,k}^{end}$. Note that this is not a certainty but in practice, if the gradient descent is making these coordinates increase, then $x_{i,k}$ will stay ``stuck'' behind $x_{j,k}$ in most cases. This suggests that the relative position structure of the point set is relatively stable during the gradient descent. This would explain why the Fibonacci set leads to the best results: it was previously observed in~\cite{PnasCDKP} that fixing the relative positions of the points using those of the Fibonacci set then optimizing leads to state-of-the-art results for the $L_{\infty}$ star discrepancy. This would provide a second argument to encourage trying to understand which structural properties make point sets have particularly low discrepancy. 

\begin{figure}[h!]
\includegraphics[scale = 0.5]{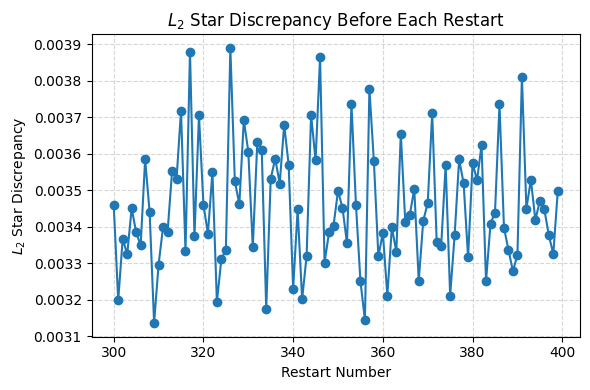}
\caption{The $L_2^*$ value immediately before restart for the last 100 restarts
of a 400 restart test initialized with a uniformly chosen random set of $n=260$. The test concluded with a best $L_2^*$ value of 0.003018 and an $L_\infty^*$ value of 0.01344 for the set corresponding to this value.}
\label{fig:restart}
\end{figure}

We also point out that to circumvent this issue we attempted to use some random restarts: once the algorithm reaches a local minimizer for gradient descent, we randomly replace some of the coordinates of the points in the set and run the algorithm again. However, this only led to very erratic behavior, and results more in line with those of separate random initializations as Figure~\ref{fig:restart} shows. It would be interesting to know how to recognize when local structure is hindering optimization, but this was outside the scope of this project.

Overall, the initial set seems to influence very strongly the quality of the output. For our experiments, we relied in dimension 2 on the Fibonacci set, while we used the Sobol' sequence or point sets from~\cite{SubsetL2} for settings where the Fibonacci set could not be used. Naturally, any practitioner who already has a high-quality set at their disposition, for example after using methods in ~\cite{SubsetL2,PnasCDKP,MPMC}, could apply the gradient descent algorithm to this point set.

\subsection{Number of iterations}

Unlike the challenges of the previous problem, the algorithm converges to a low-discrepancy set relatively quickly regardless of the initial set, number of points or discrepancy measure used.

\begin{figure}[h!]
    \centering
    \includegraphics[width=0.32\linewidth]{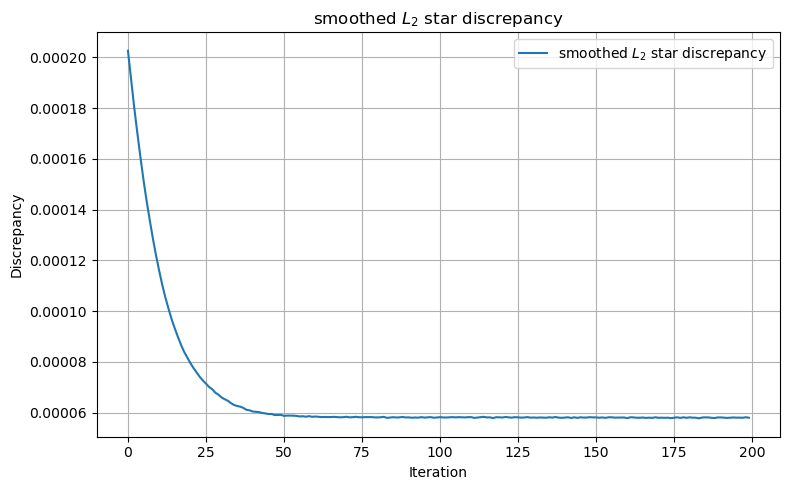}
    \hfill
    \includegraphics[width=0.32\linewidth]{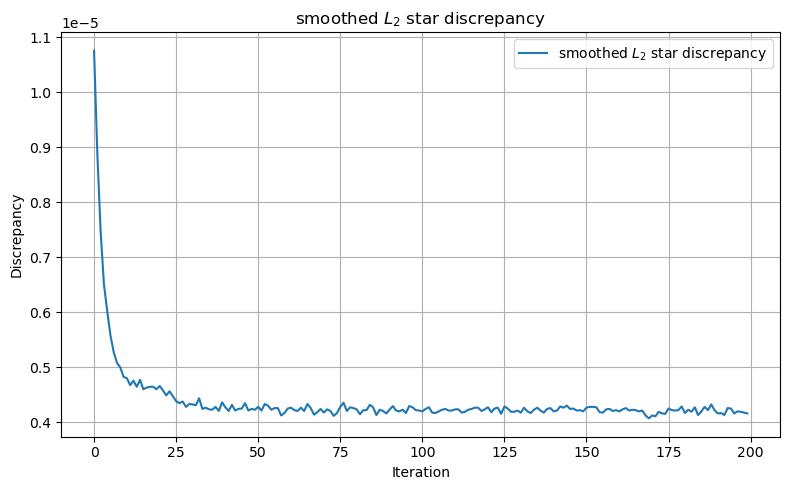}
    \hfill
    \includegraphics[width=0.32\linewidth]{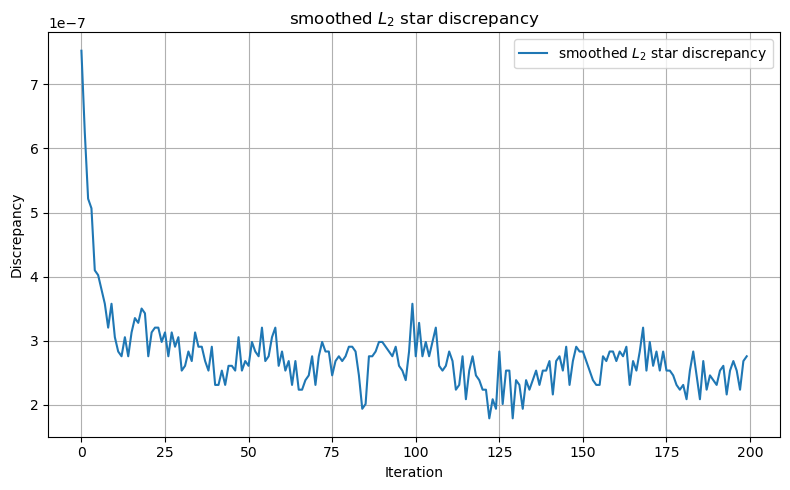}
    \caption{$L_2$ star discrepancy evolution for the Fibonacci set for $n=60,~240,~1020$ (left to right) when using the projected gradient descent algorithm. The learning rate is set to $0.0005$ for $n=60$ and $0.0001$ for the two others.}
    \label{fig:compaFibo}
\end{figure}

\begin{figure}[h!]
    \centering
    \includegraphics[width=0.5\linewidth]{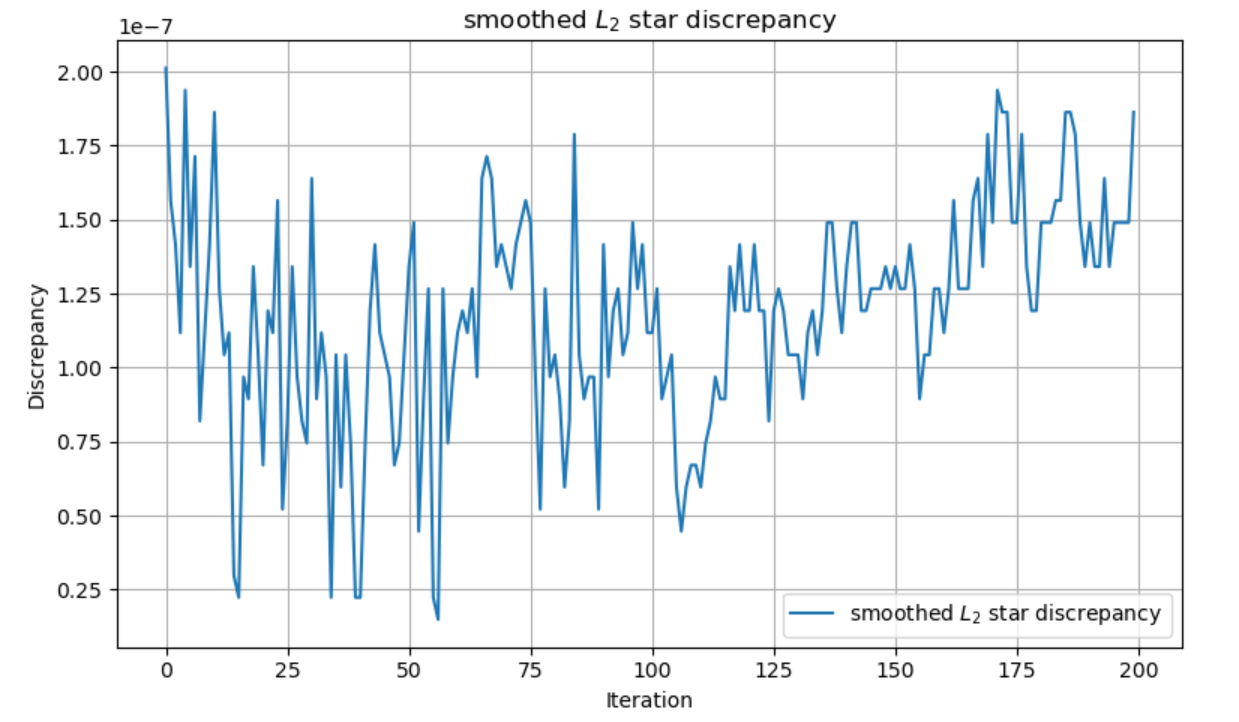}
    \caption{Gradient descent on a set of 1750 points. The fluctuations are too great to allow convergence.}
    \label{fig:fluct}
\end{figure}

Figure~\ref{fig:compaFibo} shows the evolution of the $L_2$ star discrepancy during the gradient descent with the Fibonacci set as initial point set. For all three cases, there is a sharp improvement at the start, followed by very quick stabilization for the two cases with lower $n$. However, for the $n=1020$ case, we notice that there is a quite erratic fluctuation. This appears to be caused by numerical imprecisions: as $n$ increases, these fluctuations become more and more important, to the point of limiting the algorithm's usefulness past $1500$ points. Figure~\ref{fig:fluct} illustrates this for $n=1750$, again for the Fibonacci set.

\begin{figure}[h!]
    \centering
    \includegraphics[width=0.47\linewidth]{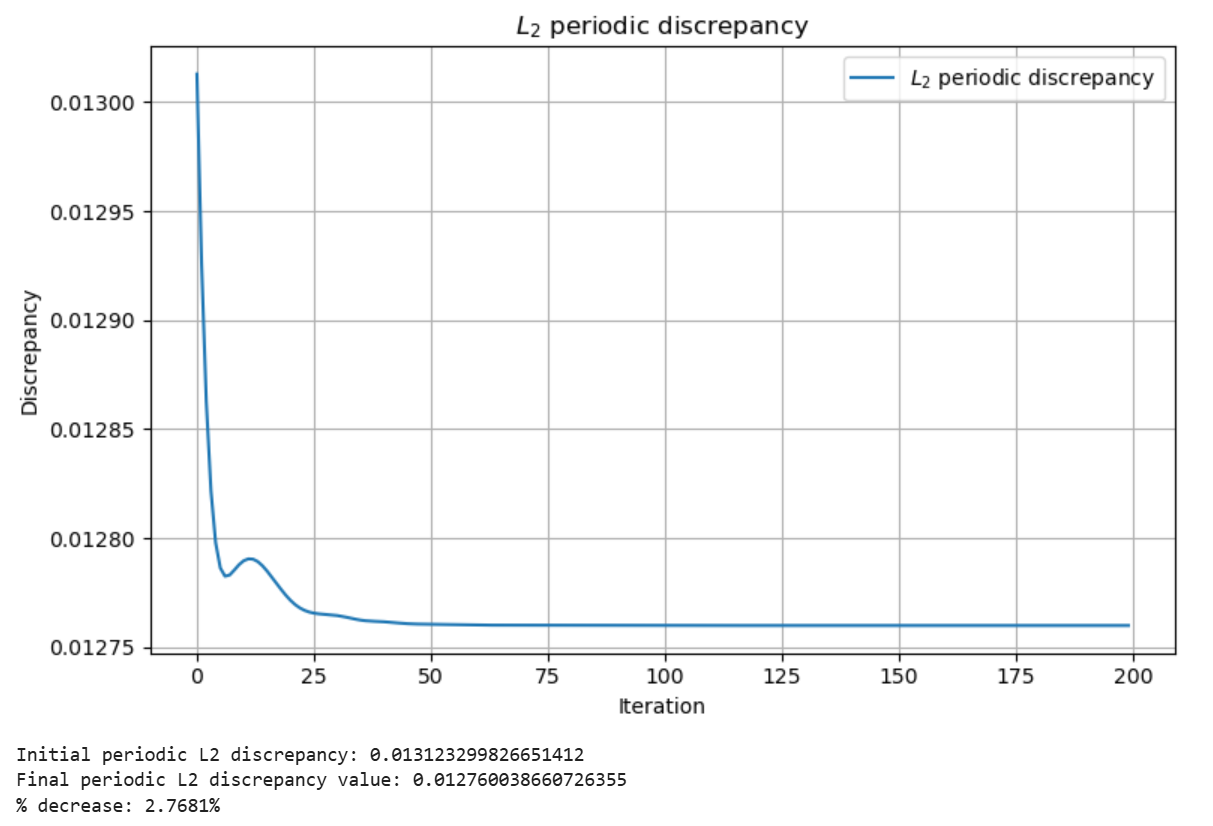}
    \hfill
    \includegraphics[width=0.47\linewidth]{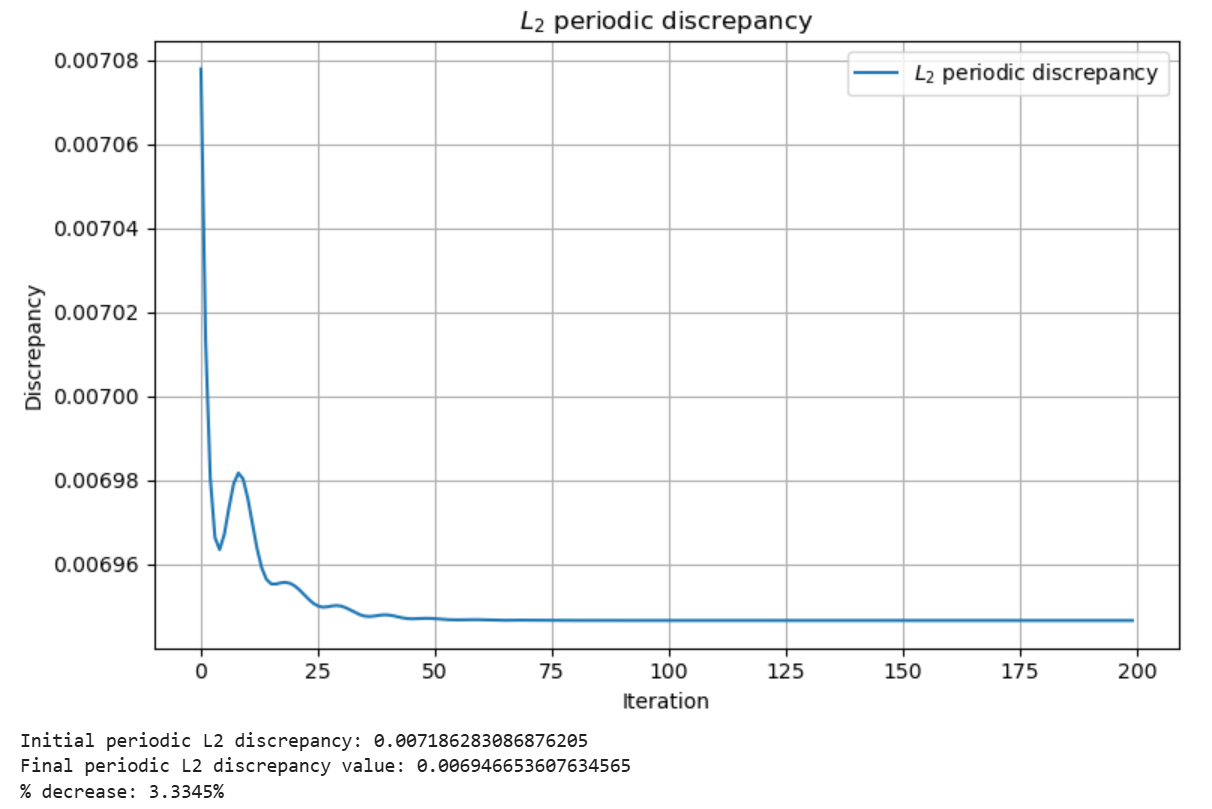}
    \caption{Periodic $L_2$ discrepancy values during the gradient descent for $n=64,~128$ points with Sobol' initialization.}
    \label{fig:perioconv}
\end{figure}

Finally, Figure~\ref{fig:perioconv} shows that this quick convergence is not limited to the $L_2$ star discrepancy. We point out that for the periodic $L_2$ discrepancy, it is well-known that integration lattices are local optimizers~\cite{Hinrichs2014}.

\section{Results}
Using the observations from the previous Section, we now run the projected gradient descent algorithm on a well-chosen initial set (Fibonacci when available, Sobol' or those of~\cite{SubsetL2} otherwise) with 200 iterations to guarantee convergence. As stated previously, each individual run took under a minute on a laptop.

\subsection{Dimension 2}

\begin{figure}[ht!]
\includegraphics[width = 0.7\linewidth]{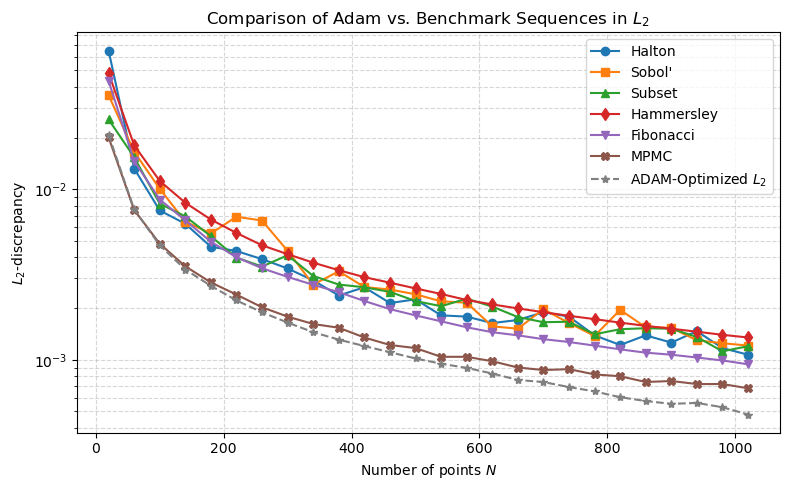}
\caption{The $L_2$ values for $n=20$ up to $n=1020$ for well-known low discrepancy sets. Projected Gradient Descent is in gray.}
\label{fig:L2disccomp}
\end{figure}

\begin{figure}[ht!]
\includegraphics[width = 0.7\linewidth]{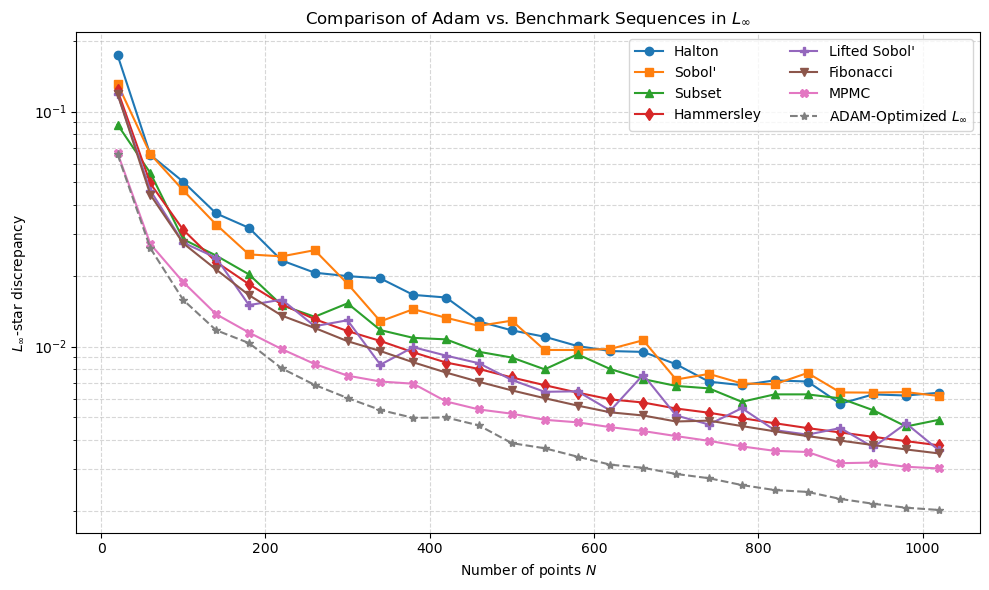}
\caption{The $L_\infty$ values for $n=20$ up to $n=1020$ for well-known low discrepancy sets. Projected Gradient Descent is in gray.}
\label{fig:Linftycompbis}
\end{figure}

We begin with the setting where we have the best benchmarks: $L_2$ and $L_{\infty}$ star discrepancy in two dimensions. We compare here our results to the discrepancies of the traditional constructions and to the results from~\cite{MPMC,PnasCDKP}. Figures~\ref{fig:L2disccomp} and~\ref{fig:Linftycompbis} respectively show the $L_2$ and $L_{\infty}$ star discrepancy values for $n=20$ to $1020$. We see that the projected gradient descent algorithm is the best performing one, in particular for the $L_{\infty}$ star discrepancy. We point out that the discrepancy values from~\cite{MPMC} were obtained with random sets for the initialization. It is highly likely that their results could be improved upon (and outperform ours) simply with a targeted initialization, as we are doing from gradient descent.

Table~\ref{tab:compaLinfty2d} then compares the results to those of~\cite{PnasCDKP}, which optimized on \emph{all} point sets whose relative positions are those of the Fibonacci set (shifted by 1, as $(0,0)$ is not a desirable point in general). While our results do not match those of~\cite{PnasCDKP}, they are competitive with them and only require a couple of minutes and a well-understood gradient descent optimization process, against the solver Gurobi~\cite{gurobi} (and longer pre-processing times for higher $n$, though finding a good set once the solver starts its heuristics is then fast).

\begin{table}[]
    \centering
    \begin{tabular}{|c|c|c|c|c|}
        \hline
        $n$& PGD (returned) & PGD (best)& MPMC&NLP\\
        \hline
        20& 0.068339 & 0.065558 & 0.0666& 0.06219  \\
        100& 0.016419 & 0.016115 &0.0188& 0.01492 \\
        180& 0.010313 & 0.010119 & 0.0115& 0.00901\\
        260& 0.007047  & 0.006965 &0.0084 & 0.00640\\
        420& 0.005072 & 0.004992 &0.0058 & 0.00412\\
        \hline
    \end{tabular}
    \caption{Comparison of the $L_{\infty}$ star discrepancy of the point set returned by the Projected Gradient Descent algorithm (PGD), the best $L_{\infty}$ found during the 200 iterations, and the results from Message-Passing Monte Carlo~\cite{MPMC} (MPMC) and non linear programming~\cite{PnasCDKP} (NLP)}
    \label{tab:compaLinfty2d}
\end{table}

\subsection{Higher dimensions}

We now turn our attention to dimensions 3, 4 and 5. We do not consider higher dimensions for two reasons. The first is that the $L_2$ discrepancy is a poorer and poorer approximation for the $L_{\infty}$ discrepancy as $d$ increases. The second is simply the lack of a good starting set to use as our initialization point. Our results already show that the Sobol' initialization does not perform particularly well in dimensions 3, 4 and 5. While we can use the sets from~\cite{SubsetL2} to circumvent this issue for $d\leq 5$, we do not have such sets in higher dimensions. As $L_2$ discrepancy comparison points are rare, and the $L_{\infty}$ star discrepancy is usually the target, we only present $L_{\infty}$ star discrepancy values in this subsection. These were computed using the algorithm from~\cite{DEM}.

Table~\ref{tab:3d} compares the results of gradient descent with those of~\cite{SubsetL2} in 3 dimensions. We can see that the Sobol' initialization performs poorly, although it is still improving on the initial Sobol' set. However, the algorithm further improves the already excellent results from~\cite{SubsetL2}. Importantly, it never makes the values noticeably worse. This confirms that the projected gradient descent algorithm can be used as a post-processing tool to complement other methods.

\begin{table}[h!]
    \centering
    \begin{tabular}{|c|c|c|c|c|}
    \hline
         $n$& PGD (Sobol')&Sobol'& PGD (L2 Subset)& L2 Subset  \\
         \hline
         50&0.088927&0.09708 & 0.05864 & 0.05952\\ 
         100&0.044201&0.06058 &0.03740 &0.03835\\
         150&0.032323&0.04483 & 0.02499&0.02612\\
         200&0.029209&0.03315 &0.02181 &0.02203\\
         250&0.020672&0.02548 & 0.01837&0.01840\\
         500&0.01194&0.01460 &0.01125&0.01207\\
         \hline
    \end{tabular}
    \caption{Comparison of $L_{\infty}$ star discrepancy values for projected gradient descent (PGD) with either Sobol' or L2 subsets from~\cite{SubsetL2} as initialization sets, with those of ~\cite{SubsetL2} and the traditional Sobol' sequence in three dimensions. All values are $L_{\infty}$ star discrepancy values.}
    \label{tab:3d}
\end{table}

Similarly, Tables~\ref{tab:4d} and~\ref{tab:5d} present the results in four and five dimensions. The observations are globally the same: the algorithm performs more poorly with a worse initial set, yet is always able to keep the quality of the set from~\cite{SubsetL2}, if not improve them outright. Note that this improvement is \emph{despite} optimizing for the $L_2$ star discrepancy, a different measure than the $L_{\infty}$ star discrepancy that is presented in the tables.
\begin{table}[h!]
    \centering
    \begin{tabular}{|c|c|c|c|c|}
    \hline
         $n$& PGD (Sobol')& Sobol'& PGD (L2 Subset)& L2 Subset  \\
         \hline
         50&0.114514&0.13422 &0.07968 & 0.08482\\ 
         100&0.068758&0.09269 &0.04660 &0.04760 \\
         150&0.054589&0.06174 &0.03937 &0.04110\\
         200&0.048119&0.05026 &0.03013 &0.03008\\
         250&0.033515&0.03822 &0.02604 &0.02596\\
         500&0.02077&0.02290 &0.01645 &0.01810  \\
         \hline
    \end{tabular}
    \caption{Comparison of $L_{\infty}$ star discrepancy values for projected gradient descent (PGD) with either Sobol' or L2 subsets from~\cite{SubsetL2} as initialization sets, with those of ~\cite{SubsetL2} and the traditional Sobol' sequence in four dimensions. All values are $L_{\infty}$ star discrepancy values.}
    \label{tab:4d}
\end{table}
\begin{table}[h!]
    \centering
    \begin{tabular}{|c|c|c|c|c|}
    \hline
         $n$& PGD (Sobol')& PGD (L2 Subset)& L2 Subset& Sobol'  \\
         \hline
         50&0.144112 &0.11355 & 0.115507&0.165488\\ 
         100&0.081881 &0.063010 &0.070071&0.120707\\
         150&0.063564 &0.05324 &0.055612&0.074899\\
         200&0.053791 &0.04249 &0.043016&0.058292\\
         500&0.02914 &0.02454 &0.026378 &0.029017\\
         \hline
    \end{tabular}
    \caption{Comparison of $L_{\infty}$ star discrepancy values for projected gradient descent (PGD) with either Sobol' or L2 subsets from~\cite{SubsetL2} as initialization sets, with those of ~\cite{SubsetL2} and the traditional Sobol' sequence in five dimensions. All values are $L_{\infty}$ star discrepancy values.}
    \label{tab:5d}
\end{table}

\subsection{Periodic and Extreme discrepancies}
Moving away from star discrepancies, we consider in this section the extreme and periodic discrepancies. Due to the lack of tailored algorithms for these measures, we only give the respective $L_2$ discrepancies. We then compare these results to the initial values of Sobol', as well as the optimized point sets in~\cite{clément2025optimizationdiscrepancymeasures} obtained via adaptations of the methods from~\cite{MPMC}.

We note that for the periodic discrepancy, the Fibonacci set introduced previously is very close to another set introduced in~\cite{BilykTY12}
$$S=\left\{\left(\frac{i}{F_k},\frac{F_{k-1}i}{F_k}\right):i \in \{0,\ldots F_{k}-1\}\right\}.$$
This is an integration lattice, and therefore a local optimum of the $L_2$ periodic discrepancy. However, it is also conjectured (and shown for $n\leq 13$ in~\cite{Hinrichs2014}) that it is the global optimizer for the $L_2$ periodic discrepancy when $n$ is a Fibonacci number. Unsurprisingly, when using this point set to initialize our algorithm, gradient descent cannot improve on the point set, though it does not move away from the optimum.

Table~\ref{tab:perio} presents the results in two dimensions for the $L_2$ periodic discrepancy. The gradient descent is able to improve both on the initial Sobol' sequence and the initial Fibonacci set, despite the latter's initial quality. As in the star discrepancy case, our method once again empirically guarantees that the discrepancy does not increase compared to that of the starting set. Finally, we note that the values indicated in~\cite{clément2025optimizationdiscrepancymeasures} for Sobol' differ from what we have found, it is possible the MPMC values do not correspond precisely to periodic discrepancy.

\begin{table}[h!]
    \centering
    \begin{tabular}{|c|c|c|c|c|c|}
    \hline
         $n$& PGD (Sobol')&PGD (Fibo)& MPMC& Sobol'&Fibo  \\
         \hline
         16 &0.048778&0.03642&0.0381&0.05163&0.03819\\ 
         32&0.02276&0.01923&0.0208 &0.02336&0.02075\\
         64&0.01276&0.01038&0.0114&0.01312&0.01158\\
         128&0.00695&0.00540&0.0060&0.00719&0.00589\\
         256&0.00409&0.00286&0.0034&0.00437&0.00317\\
         \hline
    \end{tabular}
    \caption{Comparison of $L_2$ periodic discrepancy values using projected gradient descent (PGD), with those obtained in~\cite{clément2025optimizationdiscrepancymeasures} and the initial Sobol' and Fibonacci set values.}
    \label{tab:perio}
\end{table}

Finally, we also look at the extreme discrepancy. While the Warnock-like formula was not introduced previously in this paper, it is given by Equation (8) in~\cite{clément2025optimizationdiscrepancymeasures}, and we use it with a smoothed minimum for the gradient descent. Table~\ref{tab:ext} describes our results. We see that the algorithm is always able to improve on the initial discrepancy value, even when that value is itself already very low, as in the case of the Fibonacci set.

\begin{table}[h!]
    \centering
    \begin{tabular}{|c|c|c|c|c|c|}
    \hline
        $n$ & PGD (Sobol')&PGD (Fibo)& MPMC&Sobol& Fibonacci \\
        \hline
        16 & 0.02026&0.01558 & 0.0159&0.02474 &0.01578\\
        32 &0.00984 &0.00844 &0.0088 &0.01108 &0.00854\\
        64 &0.00539 &0.00452 &0.0049 &0.00630 &0.00456\\
        128 &0.00304 &0.00240 &0.0027 &0.00347 &0.00242\\
        256 &0.00166 &0.00127 &0.0015 &0.00214 &0.00128\\
        \hline
    \end{tabular}
    \caption{$L_2$ extreme discrepancy values obtained via projected gradient descent with Sobol' and Fibonacci initial sets, compared to the initial values and those from~\cite{clément2025optimizationdiscrepancymeasures}.}
    \label{tab:ext}
\end{table}

\section{Conclusion}
We presented in this paper a projected gradient descent algorithm to optimize $L_2$ discrepancies in low dimensions. While it is highly dependent on the initial point set, our results show that the algorithm is able to substantially improve the discrepancy when possible and, if not, at least keep the quality of the starting set. It can also be applied to discrepancies with a Warnock-like formula, as we have shown for star, extreme and periodic. This method has the advantage of providing a computationally cheap method to improve a point set, in particular, it can be used as a complement of other methods such as ~\cite{SubsetL2,CDP23,MPMC}.

One of the key questions raised by this work is how to determine which local structures are desirable to construct low-discrepancy sets. Though less obvious than in~\cite{PnasCDKP}, relative position appears to play a key role in the quality of the output set. Understanding which local structures influence decisively the global property of the $L_{\infty}$ star discrepancy seems decisive in further improving low-discrepancy constructions, and this paper provides a much simpler way to test these structures than the optimization process of~\cite{PnasCDKP}.

% \bibliographystyle{alpha}
%\bibliography{refs}
\printbibliography

\end{document}